\documentclass[a4paper,11pt]{article}
\usepackage{latexsym}
\usepackage{amssymb}
\usepackage{amsfonts}
\usepackage{amsmath}
\usepackage{indentfirst}
\usepackage{graphicx}
\usepackage{epsfig}

\newtheorem{example}{Example}
\newtheorem{lemma}{Lemma}

\newtheorem{teor}{Theorem}[section]
\date{}

\author{E. Deutsch\thanks{Department of Mathematics,
Polytechnic University, Brooklyn, NY 11201, United States.} \and
E. Pergola$^{\dag}$ \and R. Pinzani\thanks{Dipartimento di Sistemi
e Informatica - Universit\`a degli Studi di Firenze, V.le G. B.
Morgagni 65, 50134 - Firenze - Italy.}}

\title{Six bijections between deco polyominoes and permutations}

\begin{document}
\maketitle \hspace{3.7cm}{To the memory of Alberto Del Lungo
(1965--2003)}
\vspace{1cm}
\begin{abstract}
In this paper we establish six bijections between a particular
class of polyominoes, called deco polyominoes, enumerated
according to their directed height by n!, and permutations. Each
of these bijections allows us to establish different
correspondences between classical statistics on deco polyominoes
and on permutations.
\end{abstract}

\section{Introduction}
Deco polyominoes were introduced in \cite{BDP} with the aim to
find a class of polyominoes counted by the factorial. Indeed, in
the formula giving the height of directed column-convex
polyominoes, see \cite{BBDP}, there is a factorial and, as the
authors say in their paper, X. Viennot called their attention to
the problem discussed in that paper. Moreover, deco polyominoes
have provided a nice tool for the random generation of
permutations.

On the other hand, permutations are a widely studied combinatorial
objects but, nevertheless, there are few bijections between them
and other combinatorial objects.

This study brings the two structures together by establishing six
different bijections between them and by finding correspondences
between certain statistics on deco polyominoes and on
permutations.

In Section \ref{sec_notations} we recall some definitions and
properties regarding deco polyominoes and permutations. Sections
\ref{sec_bij_1}, \ref{sec_bij_2}, \ref{sec_bij_3},
\ref{sec_bij_4}, \ref{sec_bij_5}, and \ref{sec_bij_6} describe and
exemplify the six bijection we establish between deco polyominoes
and permutations, pointing out correspondences between certain
statistics defined on them.

\section{Notations and Definitions}\label{sec_notations}
We introduce some notations and definitions regarding deco
polyominoes and permutation, respectively.
\subsection{Deco polyominoes}\label{sec deco polyominoes}
\begin{figure}[!hbp]
\begin{center}
\includegraphics[scale=.385]{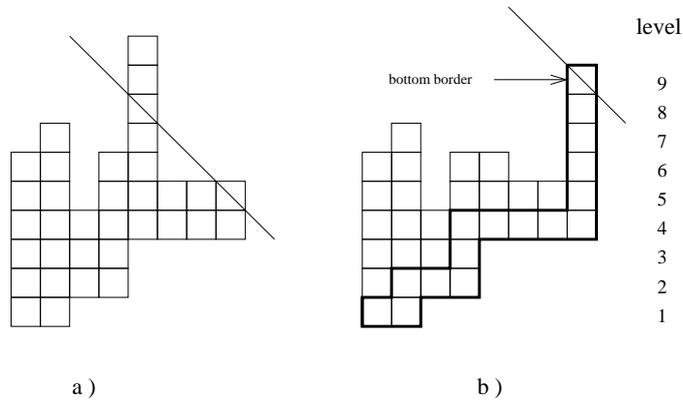}
\end{center}
\caption{a) A directed column-convex polyomino. b) A polyomino of
height $16$, width $8$, area $34$, whose last column level is $9$}
\label{polyomino}
\end{figure}

In the $\Re^2$ plane, a \emph{cell} is a unitary square
$[i,i+1]\times[j,j+1]$, $i,j \in \aleph$, and a \emph{polyomino}
is a connected set of pairs of cells having one side in common.
Polyominoes are defined up to a  translation. We can obtain a
\emph{directed} polyomino by starting out from a cell, called
\emph{source}, and by adding other cells in predetermined
directions, such as East and North, that is, to the right of or
over existing cells. In this way, a polyomino grows in a
\emph{privileged direction}. A \emph{column} (resp. \emph{row)} is
the intersection of a polyomino with an infinite vertical (resp.
horizontal) strip $[i,i+1]\times\Re$ (resp. $\Re\times [j,j+1]$).
A \emph{directed column-convex} polyomino is a directed polyomino
whose columns are connected (see Figure \ref{polyomino} a)).
Finally, the \emph{directed height} of a directed polyomino is the
number of lines orthogonal to the privileged direction that go
through the cell center (from here on we call it simply the \emph{
height}), while its \emph{vertical height} is the number of rows
and its \emph{width} is the number of columns. The \emph{area} of
a polyomino is defined to be the number of its cells. The
\emph{level} of a column of a directed polyomino is the number of
rows lying between the top of the column itself and the bottom of
the source. The \emph{bottom border} of a directed polyomino is
made up of the cells lying on the lowest path going from the
source to the highest rightmost cell (see Figure \ref{polyomino}
b)). A \emph{parallelogram polyomino} can be defined as an array
of unit cells that is bounded by two lattice paths which use the
steps $(1,0)$ and $(0,1)$, and which intersect only at their
origin and extremity. They form a subclass of convex (i.e. row-
and column-convex) polyominoes. It is well-known that
parallelogram polyominoes of semiperimeter $n+1$ are counted by
the Catalan number $\frac{1}{n+1}{\left( \begin{array}{c} 2n\\ n\\
\end{array} \right)}$, see \cite{DV} for example. We will consider a particular class of
directed column-convex polyominoes, called \emph{deco polyominoes}
after the French \emph{derni\`ere colonne}: last column. They are
defined as directed column-convex polyominoes in which the height
is attained only in the last column (see Figure \ref{polyomino}
b)). Deco polyominoes have been introduced by E. Barcucci, A. Del
Lungo, and R. Pinzani \cite{BDP} (see also \cite{BBD}). We denote
the class of deco polyominoes of height $n$ by $D_n$. Clearly, the
height of a deco polyomino is equal to the number of cells in its
bottom border. Also, for any deco polyomino of height $n$ we have
\emph{width} $+$ \emph{level of last column} $=$ $n+1$. It is easy
to see that the parallelogram polyominoes of semiperimeter $n+1$
form a subset of the deco polyominoes of height $n$.

\begin{figure}[!hbp]
\begin{center}
\includegraphics[scale=.385]{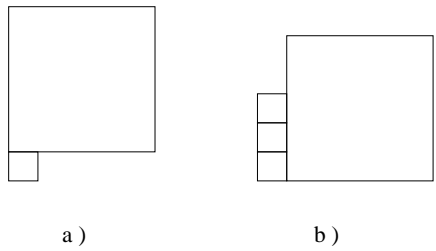}
\end{center}
\caption{Decomposition of the deco polyominoes}
\label{construction_polyomino}
\end{figure}

As shown in \cite{BDP}, the set $D_n$ of deco polyominoes of
height $n$ ($n\geq 2$) admits the following decomposition. If a
deco polyomino of height $n$ has no cell to the right of the
source, then it is made up of the source and a deco polyomino of
height $n-1$ attached to the North side of the source (see Figure
\ref{construction_polyomino} a)). Otherwise, when there is a cell
to the right of the source, the polyomino is made up of the first
column, containing $k$  cells , where $k \leq n-1$ (because height
$n$ is attained only in the last column), and a deco polyomino of
height $n-1$, attached to the East side of the source (see Figure
\ref{construction_polyomino} b)).

From the above decomposition of $D_n$, taking into account that in
the first case we have $|D_{n-1}|$ possibilities and in the second
case we have $(n-1)|D_{n-1}|$ possibilities, we obtain $|D_n| =
n|D_{n-1}|$; this, together with $|D_1|=1$, yields $|D_n| = n!$.

A deco polyomino of height $n$ can be built step by step by a
sequence  of $n$ steps of two kinds:

\emph{elevation}: add a cell at the bottom of the leftmost column
of the previous deco polyomino; the first step is always of this
kind, the "previous" deco polyomino being the empty one;

\emph{column pasting}: add a new column to the left of the
previous deco polyomino in such a way that the bottoms of the
first two columns lie at the same level.

\begin{figure}[!hbp]
\begin{center}
\includegraphics[scale=.385]{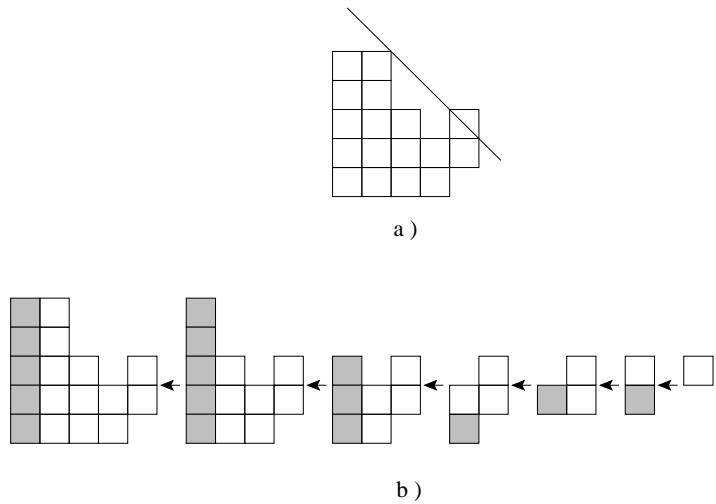}
\end{center}
\caption{a) A deco polyomino. b) The step-by-step construction of
the polyomino in a)} \label{BDP construction}
\end{figure}

For a deco polyomino of height $n$ we define $a_j = 0$ if the
$j$-th step of the above described construction is an elevation
and $a_j = k$ if it is a pasting of a column of length $k$. Thus,
a deco polyomino of height $n$ can be coded by the sequence $(a_n,
a_{n-1},\ldots, a_2, a_1)$ (it is convenient to list the $a_j$'s
in this order). Since after step $j$ of the step-by-step
construction we obtain a deco polyomino of height $j$, it follows
that $0 \leq a_j \leq j-1$ for $j=1,\ldots,n$. For example, for
the deco polyomino of Figure \ref{BDP construction} a), the
step-by-step construction is shown in Figure \ref{BDP
construction} b) (from right to left) and, consequently, the
corresponding code is $(5,5,3,0,1,0,0)$.

A different coding of a deco polyomino is given in \cite{BDP}; the
reader may derive the simple connection between the two codings.

\subsection{Permutations}\label{def_permutation}
We denote the set of all permutations of the set
$[n]=\{1,2,\ldots,n\}$ by $S_n$.

We say that a permutation $\pi \in S_n$ \emph{contains a
subsequence of type} $\tau \in S_k$ if there exists a sequence of
indices $1 \leq i_1 \le i_2 \le \ldots \le i_k\leq n$ such that
$\pi(i_1)\pi(i_2)\ldots \pi(i_k)$ is ordered as $\tau$. If no such
sequence exists, then the permutation $\pi$ is said
$\tau$-avoiding. It is well known that $321$-avoiding permutations
are counted by Catalan numbers (see \cite{B}, p. 133).

\begin{example}
The permutation $351264$ is $321$-avoiding.
\end{example}

By a \emph{sequence} of length $n$ we mean a list of $n$ distinct
positive integers. To each sequence $s$ of length $n$ we associate
a permutation of $[n]$ in a natural way: we relabel the smallest
number in $s$ as $1$, the second smallest number as $2$, and so
on, relabelling the largest number in $s$ as $n$. We call this
permutation the \emph{reduction} of the sequence $s$, denoted
$red(s)$.

\begin{example}
Let $s=572396$ then $red(s)=351264$.
\end{example}

For a permutation $\pi=\pi_1\pi_2\ldots \pi_n$, we define the
\emph{reverse} of $\pi$ as the permutation
$\pi^r=\pi_n\pi_{n-1}\ldots \pi_1$ and the \emph{complement} of
$\pi$ as the permutation having entries $n+1-\pi_i$,
$i=1,\ldots,n$.

Let $\pi=\pi_1\pi_2\ldots \pi_n$ be a permutation. We say that $i$
is a \emph{descent} of $\pi$ if $\pi_i>\pi_{i+1}$. If $\pi$ has
$k-1$ descents, then $\pi$ is the union of $k$ increasing
subsequences of consecutive entries. These are called the
\emph{ascending runs} of $\pi$. In the same way one defines the
concepts of \emph{ascent} and \emph{descending runs}. An entry in
a permutation which is smaller than all the entries that follow it
is called \emph{right-to-left minimum} (see \cite{B}, p. 98).
Clearly, the right-to-left minima form an increasing sequence.

\begin{example}
The permutation $2371546$ has descents $3, 5$, ascents $1, 2, 4,
6$, three ascending runs $237, 15, 46$, five descending runs $2,
3, 71, 54, 6$ and right-to-left minima $6, 4, 1$.
\end{example}

\bigskip

In a permutation $\pi=\pi_1\pi_2\ldots \pi_n$, an \emph{inversion}
is a pair $i<j$ such that $\pi_i>\pi_j$. The number of inversions
of $\pi$ will be denoted by $inv(\pi)$.

If $c_i$ is the number of $j>i$ with $\pi_j<\pi_i$, then
$(c_1,c_2,\ldots,c_n)$ is called the \emph{right inversion vector}
of $\pi$ (\cite{J} p. 188, where the terms inversion vector and
inversion table are used). Clearly, $0\leq c_i\leq n-i$. The right
inversion vector is sometimes called code  (\cite{BJS} p. 357;
\cite{M} p. 9).

Clearly, $inv(\pi)=\Sigma_{i=1}^n c_i$.

It is known that the right inversion vector determines uniquely
the permutation.

\begin{example}
If $\pi=53728146$, then the right inversion vector is
$(4,2,4,1,3,0,$ $0,0)$, and $inv(\pi)=14$.
\end{example}

Carlitz \cite{C} (see also \cite{S}) defines the statistics
$\emph{inv}_c$ on $S_n$ as follows: express $\pi\in S_n$ in
standard cycle form (i.e. cycles ordered by increasing smallest
elements in the first position); then remove the parentheses and
count the inversions in the obtained permutation.

\bigskip

\begin{example}
If $\pi=2357146=(1235)(476)$, then $inv_c(\pi)=2$, the number of
inversions in the permutation $1235476$.
\end{example}

Clearly, the mapping defined above, which associates to any
permutation the permutation obtained by removing the parentheses,
is not a bijection (all images start with the entry $1$).

\section{Bijection No. 1}\label{sec_bij_1}

This bijection, say $\Phi_1$, has been introduced in \cite{BDP}.
It is defined recursively in the following manner.

\begin{figure}[!hbp]
\begin{center}
\includegraphics[scale=.385]{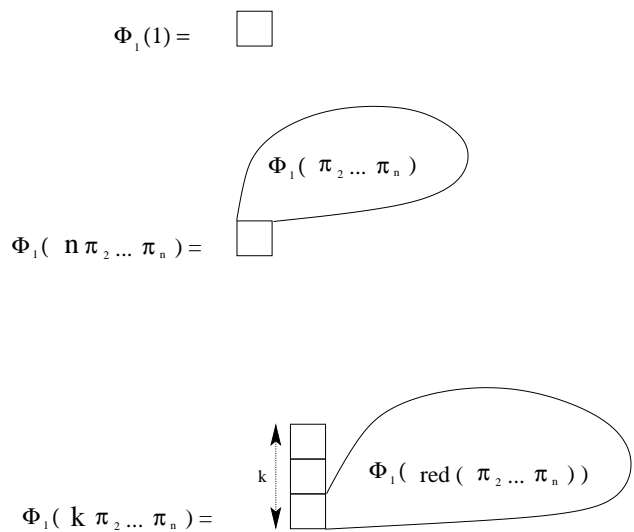}
\end{center}
\caption{A graphical representation of bijection $\Phi_1$}
\label{bij_1}
\end{figure}

To the permutation $1$ of $S_1$ there corresponds the single
cycle-cell polyomino. Let $\pi=\pi_1\pi_2\ldots \pi_n$, $n \geq
2$. If $\pi_1=n$, then the image of $\pi$ is defined to be the
polyomino obtained by attaching the polyomino corresponding to
$\pi_2\pi_3\ldots \pi_n$ to the North side of the single-cell
polyomino. If $\pi_1=k<n$, then the image of $\pi$ is defined to
be the polyomino obtained by attaching the polyomino corresponding
to the permutation $red(\pi_2,\pi_3,\ldots,\pi_n)$ to the East
side of a column of $k$ cells (see Figure \ref{bij_1}).

It is straightforward to describe the inverse mapping of this
bijection.

\begin{figure}[!hbp]
\begin{center}
\includegraphics[scale=.385]{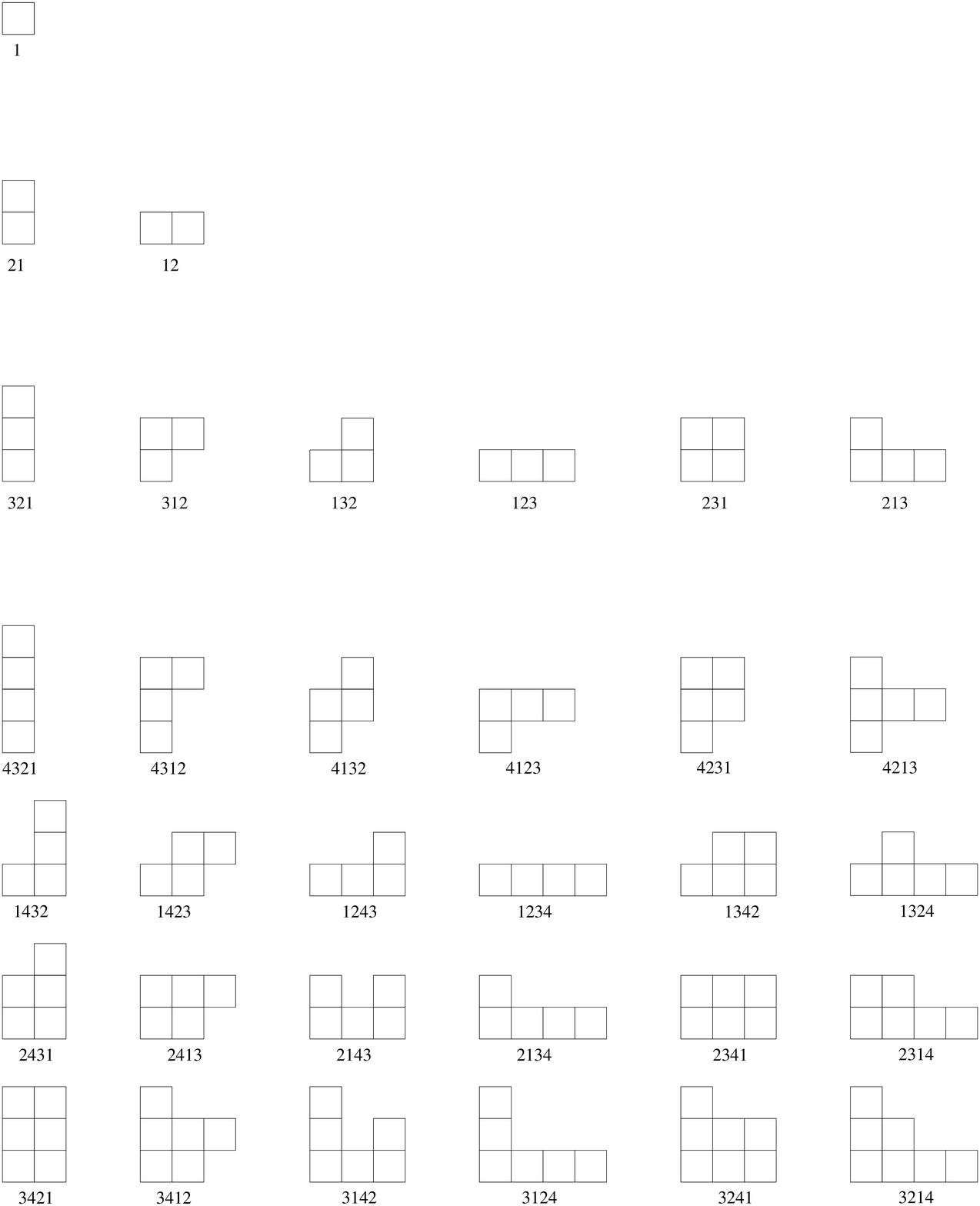}
\end{center}
\caption{Bijection $\Phi_1$ for $n=1,2,3,4$} \label{ex bij_1}
\end{figure}

Figure \ref{ex bij_1} shows this bijection for $n=1, 2, 3, 4$.

\bigskip

From the definition of this bijection it follows that the number
of cells in the last column of the deco polyomino $\Phi_1(\pi)$ is
equal to the length of the first ascending run of $\pi^r$, the
reverse of the corresponding permutation.

\section{Bijection No. 2}\label{sec_bij_2}

We present here a simple bijection, say $\Phi_2$, between
permutations in $S_n$ and deco polyominoes of height $n$, based on
the coding we have derived from the step by step construction of a
deco polyomino at the end of Section \ref{sec deco polyominoes}.

\begin{figure}[!hbp]
\begin{center}
\includegraphics[scale=.4]{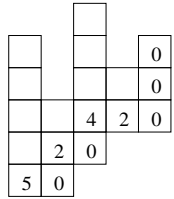}
\end{center}
\caption{A deco polyomino and its code according to bijection
$\Phi_2$} \label{deco}
\end{figure}

We start with $\Phi_2^{-1} : D_n \rightarrow S_n$. Let $\delta$ be
a given deco polyomino  and let $(b_n, b_{n-1},\ldots, b_2, b_1)$
be its code. It is convenient to insert these numbers in the
appropriate cells: the $0$'s in the new cells that produce the
elevation and the $k$'s ($k> 0$) in the bottom cells of the pasted
columns (see Figure \ref{deco}).

Since, as we have seen, $0 \leq b_j \leq j-1$, if we denote $c_j =
b_{n+1-j}$, then $0 \leq c_j \leq n-j$. Consequently, $b_n,
b_{n-1},\ldots, b_2, b_1$ $(=c_1, c_2,\ldots, c_n)$ can be viewed
as a right inversion vector. To it there corresponds a unique
permutation $\pi = \Phi_2^{-1}(\delta)$.

\begin{example}
The deco polyomino of Figure \ref{deco} is coded by the sequence
$(5,0,2,$ $0,$$4,2,0,0,0)$ and to this right inversion vector
there corresponds the permutation $614297358$.
\end{example}

It is straightforward to obtain the deco polyomino  $\Phi_2(\pi)$
from a given permutation $\pi$: find the right inversion vector
$(c_1, c_2,\ldots, c_n)$ of $\pi$ and, viewing this as a code,
perform the step-by-step construction of the corresponding deco
polyomino.

\begin{figure}[!hbp]
\begin{center}
\includegraphics[scale=.4]{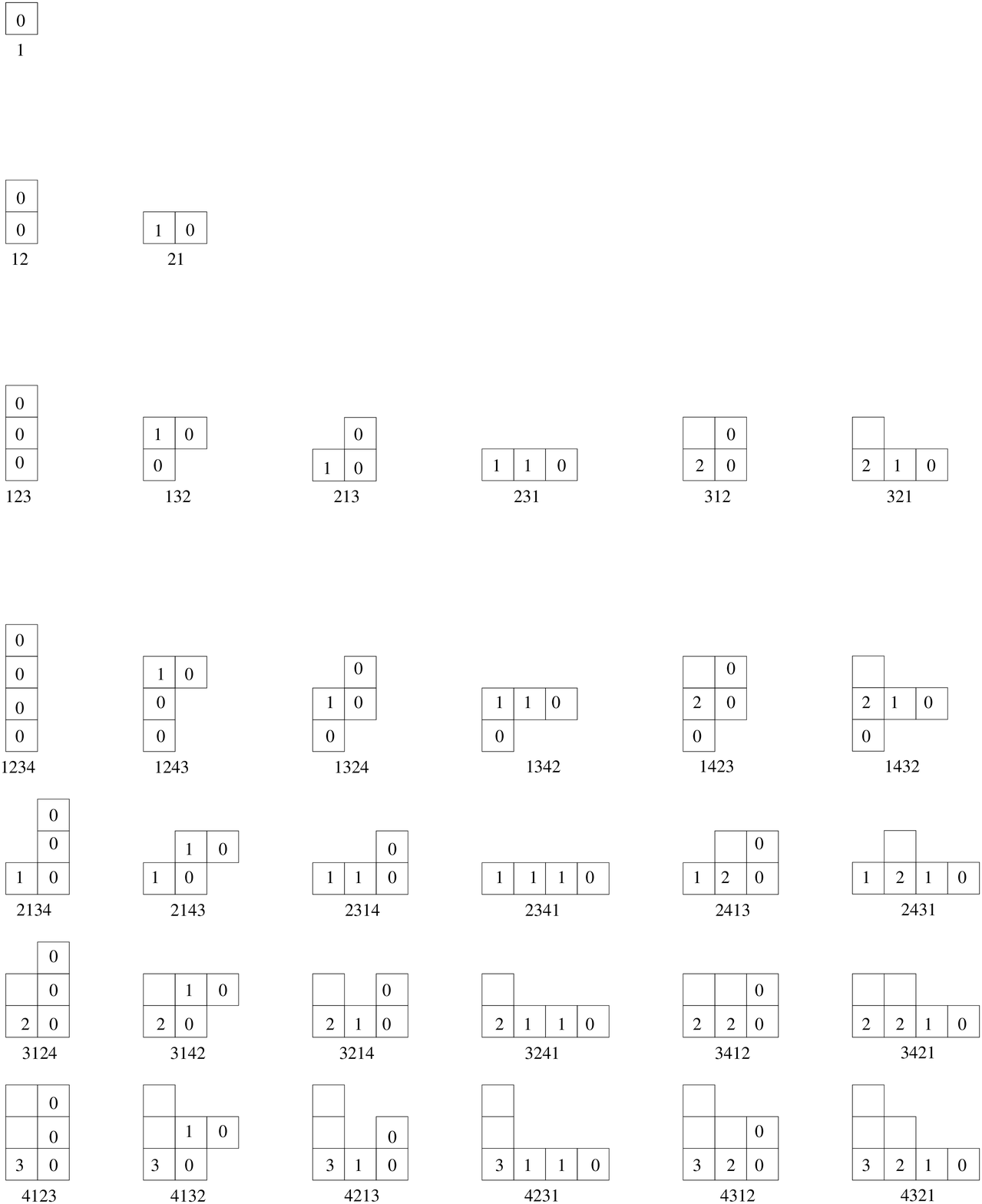}
\end{center}
\caption{Bijection $\Phi_2$ for $n = 1, 2, 3, 4$} \label{ex bij_2}
\end{figure}

Figure \ref{ex bij_2} shows this bijection for $n=1, 2, 3, 4$.

The following relations between a permutation $\pi \in S_n$ and
its corresponding deco polyomino $\Phi_2(\pi) \in D_n$ are
immediate:
\begin{itemize}
    \item the level of the last column of $\Phi_2(\pi)$ is equal to the number of
    right-to-left minima of $\pi$;
    \item the number of cells in the last column of $\Phi_2(\pi)$ is
    equal to the length of the last ascending run of $\pi$;
    \item the area of $\Phi_2(\pi)$ is equal to  $inv(\pi) +$ number of right-to- left minima of
    $\pi$;
    \item if $m_1,m_2,\ldots,m_r$ are the positions of the
    right-to-left minima of $\pi$, then the lengths of the rows of
    the bottom border of $\Phi_2(\pi)$, starting from the lowest one,
    are $m_1,m_2-m_1,\ldots,m_r-m_{r-1}$, respectively.
\end{itemize}

We also note that to the identity permutation in $S_n$ there
corresponds the deco polyomino that consists of a single column of
$n$ cells; clearly, the step-by-step construction of this
polyomino involves only elevations. On the other hand, for an
arbitrary permutation  $\pi \in S_n$,  $inv(\pi)$  is equal to the
number of those cells in   $\Phi _2(\pi)$  that have been added by
column pasting.

\begin{teor}
A permutation $\pi$ is $321$-avoiding if and only if $\Phi_2(\pi)$
is a parallelogram polyomino.
\end{teor}
\textbf{Proof:} Taking into account the interpretation of the
$b_k$'s in the step-by-step construction of a deco polyomino, it
follows that we obtain a parallelogram polyomino if and only if
the number of $0$'s between two consecutive nonzero $b_i$ and
$b_j$ is at least $b_i-b_j$ (otherwise, the top of the column
corresponding to $b_i$ is higher than the top of the column
corresponding to $b_j$). But this is exactly the necessary and
sufficient condition for a permutation to be $321$-avoiding, given
by Billey, Jockush, and Stanley (\cite{BJS}, Theorem 2.1).$\hfill
\square$

\section{Bijection No. 3}\label{sec_bij_3}

We now define a new bijection, say $\Phi_3: S_n \rightarrow D_n$,
recursively. We describe $\Phi_3^{-1}$.

\begin{figure}[!hbp]
\begin{center}
\includegraphics[scale=.4]{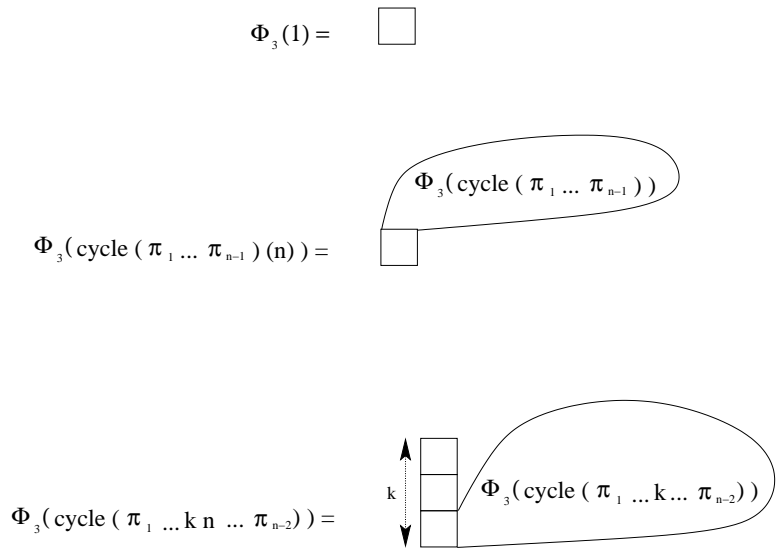}
\end{center}
\caption{A graphical representation of bijection $\Phi_3$}
\label{bij_3}
\end{figure}

To the deco polyomino of height $1$ there corresponds the
permutation $1$. Let $\delta \in D_n$. If $\delta$ is obtained by
elevation from $\delta' \in D_{n-1}$, then $\Phi_3^{-1}(\delta)$
is obtained by adding $n$ as a new cycle to the cycle form of
$\Phi_3^{-1}(\delta')$. If $\delta$ is obtained by column pasting
from $\delta' \in D_{n-1}$, then $\Phi_3^{-1}(\delta)$ is obtained
by inserting $n$ into the cycle form of $\Phi_3^{-1}(\delta')$ on
the immediate right of $k$, where $k$ is the length of the pasted
column (see Figure \ref{bij_3}).

Again, it is straightforward to describe the inverse mapping.

\begin{figure}[!hbp]
\begin{center}
\includegraphics[scale=.4]{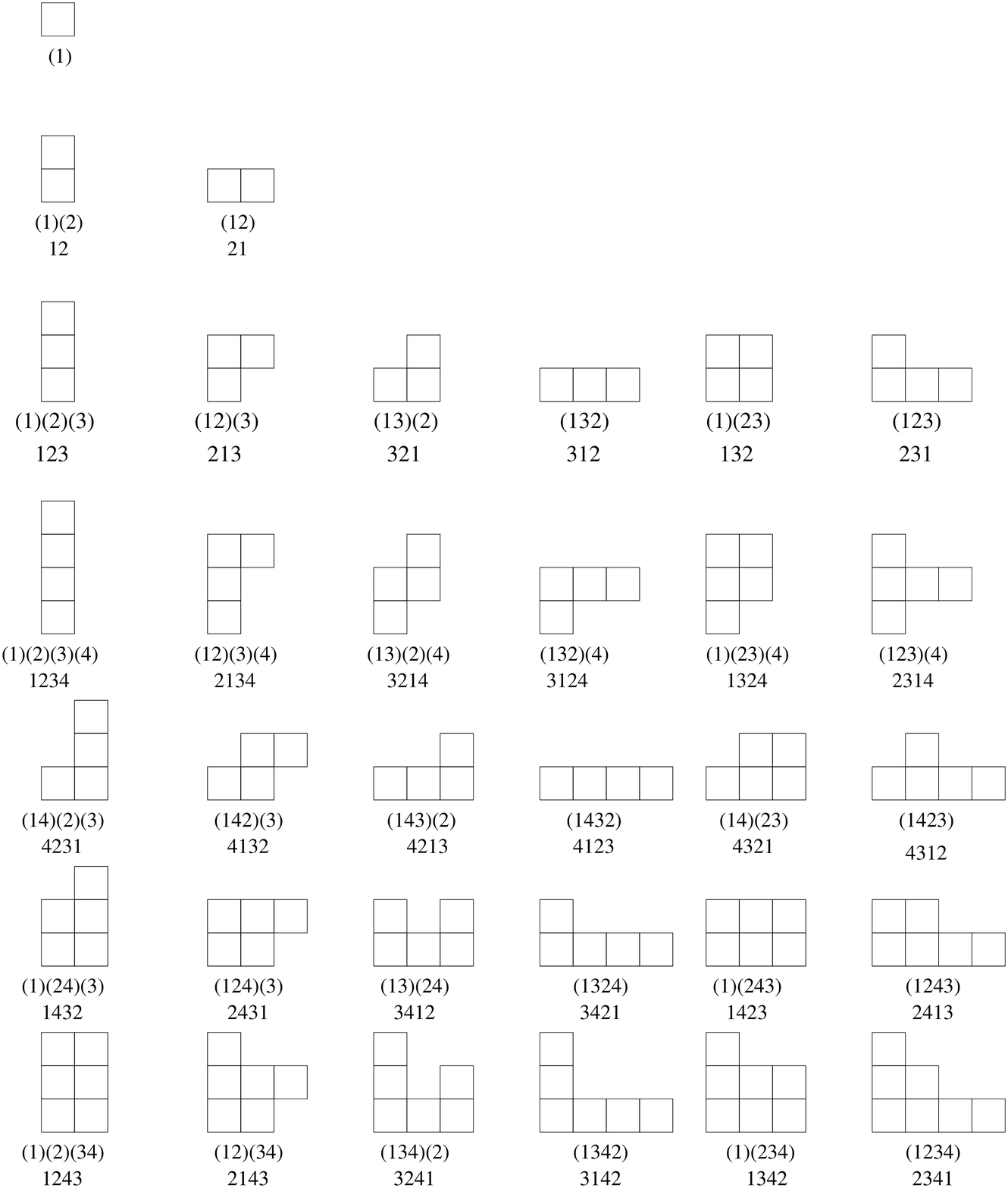}
\end{center}
\caption{Bijection $\Phi_3$ for $n=1, 2, 3, 4$} \label{ex bij_3}
\end{figure}

Figure \ref{ex bij_3} shows this bijection for $n=1, 2, 3, 4$.

It is easy to see that the number of cycles of a permutation $\pi$
is equal to the level of the last column of $\Phi_3(\pi)$. Indeed,
this is true for the permutation $1 \in S_1$ and new cycles are
obtained only through elevation. Equivalently, the width of a deco
polyomino of height $n$ is equal to $n+1-s$, where $s$ is the
number of cycles of the corresponding permutation.

\section{Bijection No. 4}\label{sec_bij_4}

\begin{figure}[!hbp]
\begin{center}
\includegraphics[scale=.4]{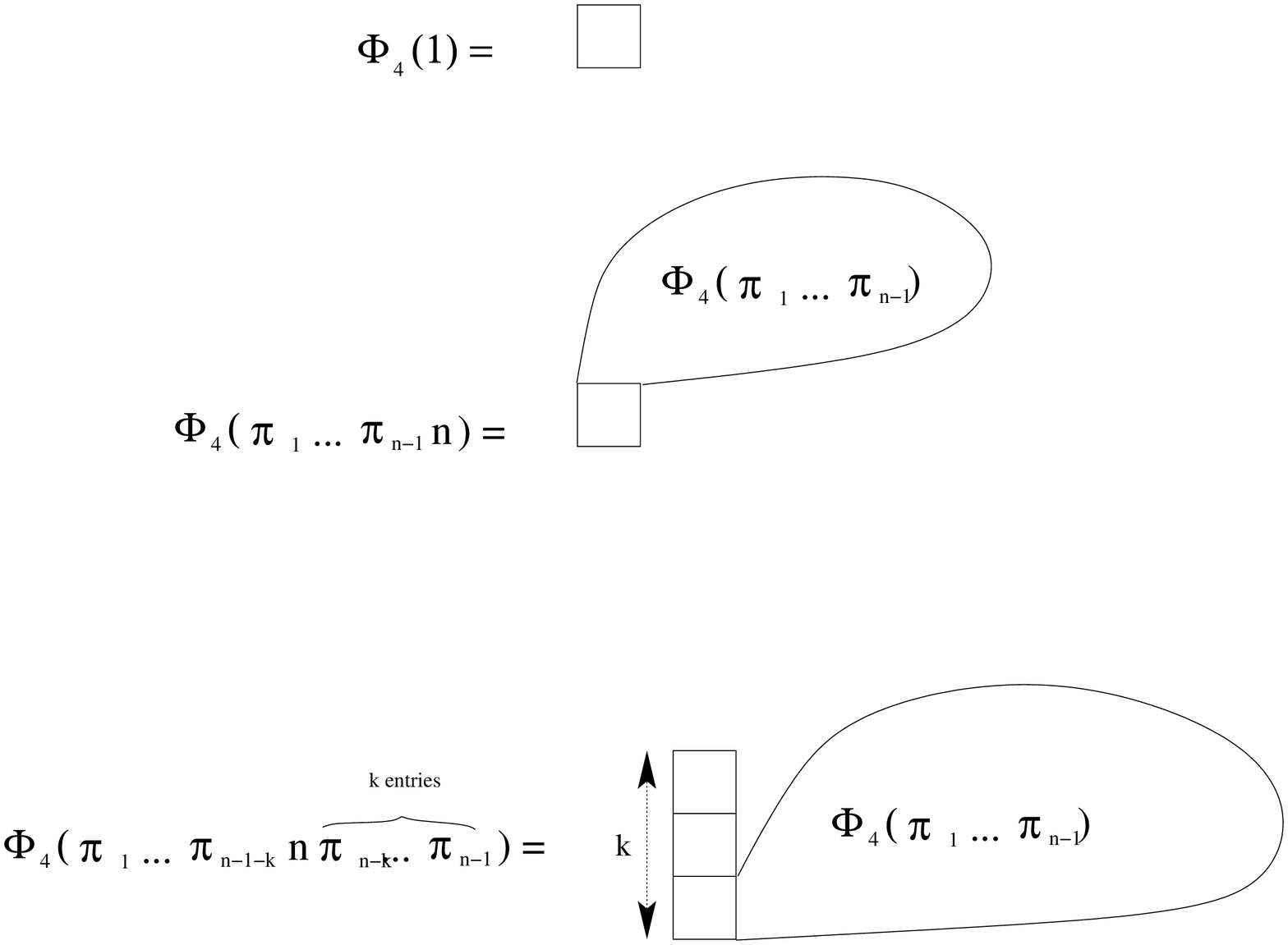}
\end{center}
\caption{A graphical representation of bijection $\Phi_4$}
\label{bij_4}
\end{figure}

We now define a new bijection $\Phi_4: S_n \rightarrow D_n$,
recursively. We describe $\Phi_4^{-1}$. To the deco polyomino of
height $1$ there corresponds  the permutation $1$. Let $\delta \in
D_n$. If $\delta$ is obtained by elevation from $\delta' \in
D_{n-1}$ and $\Phi_4^{-1}(\delta')=\pi_1 \pi_2 \ldots  \pi_{n-1}$,
then we define $\Phi_4^{-1}(\delta) =\pi_1 \pi_2 \ldots \pi_{n-1}
n$. If $\delta$ is obtained by column pasting from $\delta' \in
D_{n-1}$, then we define $\Phi_4^{-1}(\delta) =\pi_1 \pi_2 \ldots
\pi_{n-1-k} n \pi_{n-k} \ldots \pi_{n-1}$, where $k$ is the length
of the pasted column (see Figure \ref{bij_4}).

Again, it is straightforward to describe the inverse mapping.

\begin{figure}[!hbp]
\begin{center}
\includegraphics[scale=.37]{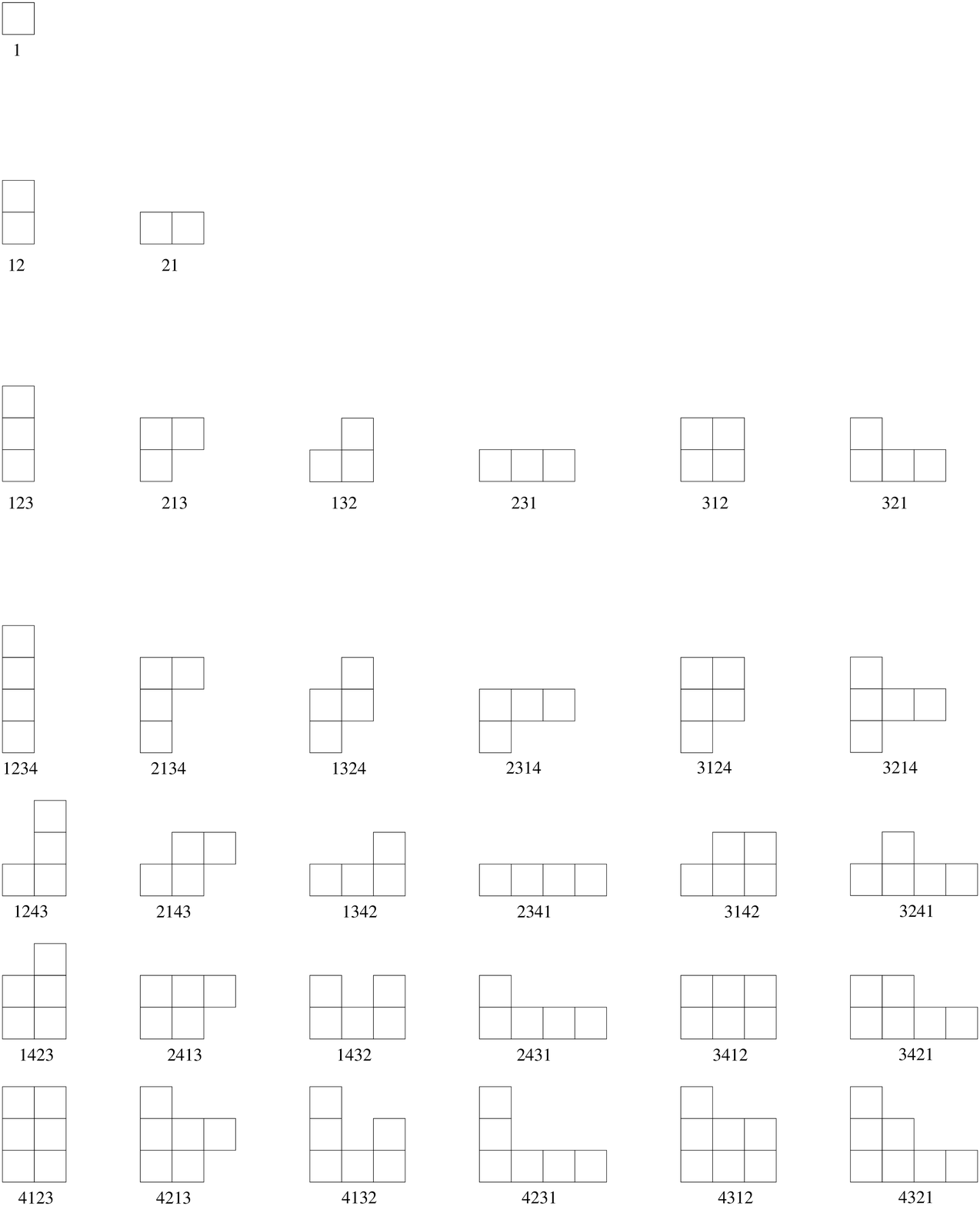}
\end{center}
\caption{Bijection $\Phi_4$ for $n=1, 2, 3, 4$} \label{ex bij_4}
\end{figure}

Figure \ref{ex bij_4} shows this bijection for $n=1, 2, 3, 4$.

It is easy to show by induction that for a permutation $\pi \in
S_n$ and its corresponding deco polyomino $\Phi_4(\pi) \in D_n$ we
have:
\begin{itemize}
\item $inv(\pi) = area(\Phi_4(\pi))-$ level of the last column of
$\Phi_4(\pi)$;

equivalently,

\item $inv(\pi) = area(\Phi_4(\pi))$ $+$ $width(\Phi_4(\pi))-$
$(n+1)$.
\end{itemize}

The proof of the next lemma, needed for the subsequent theorem, is
left to the reader.

\begin{lemma}
Let $\delta$ be a parallelogram polyomino and let $\pi$ be the
corresponding permutation ($\pi = \Phi_4^{-1}(\delta)$). Then the
length of the first column of $\delta$ is equal to the length of
the last ascending run of $\pi$.
\end{lemma}

\begin{teor}
A permutation $\pi$ is $321$-avoiding if and only if $\Phi_4(\pi)$
is a parallelogram polyomino.
\end{teor}
\textbf{Poof:} Let $\delta$ be a deco polyomino of height $n$, let
$\delta_1, \delta_2,\ldots,\delta_n=\delta$ be the deco
polyominoes obtained successively by the step-by-step construction
of $\delta$, and let $\pi_i$ be the permutation corresponding to
$\delta_i$, $i=1,\ldots,n$ under the considered bijection. As long
as $\delta_i$ consists of a single column, the permutation $\pi_i$
is the identity permutation on $S_i$. As long as $\delta_i$
consists of two columns, $\delta_i$ is a parallelogram polyomino
(because height is attained only in the last column) and $\pi_i$
is $321$-avoiding because it has exactly one descent. In the
subsequent steps, if any, we obtain for the first time a deco
polyomino $\delta_j$ that is not a parallelogram polyomino only
following a pasting of a column of length greater than the first
column of $\delta_{j-1}$. But, due to the lemma, this is the only
possibility for the $321$-avoiding permutation $\pi_{j-1}$ to go
into a permutation $\pi_j$ containing the pattern $321$. $\hfill
\square$

\section{Bijection No. 5}\label{sec_bij_5}

\begin{figure}[!hbp]
\begin{center}
\includegraphics[scale=.4]{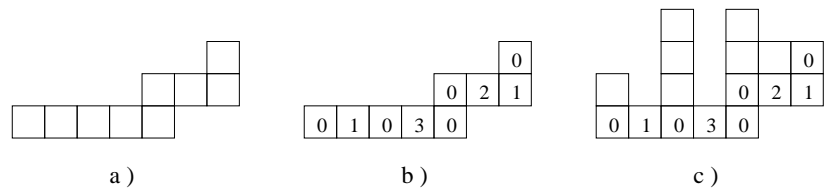}
\end{center}
\caption{An example for bijection $\Phi_5$} \label{bij_5}
\end{figure}

Parallel with the presentation of the next bijection, say
$\Phi_5:S_n\rightarrow D_n$, for an arbitrary permutation $\pi \in
S_n$, we show its steps on the example $\pi=372196458$. We write
$\pi$ in standard cycle form, $\pi=(13274)(598)(6)$. The length of
these cycles ($5$, $3$, and $1$) will be the lengths of the rows
of the bottom border of the corresponding deco polyomino (see
Figure \ref{bij_5} a) ).
Removing the parentheses in the last expression for $\pi$, we
obtain the auxiliary permutation $\pi'=132745986$. We consider its
right inversion vector $(0,1,0,3,0,0,2,1,0)$. Clearly, the sum of
these numbers (i.e. the inversion number of $\pi'$) is the Carlitz
inversion number of $inv_c(\pi)$, defined in Section
\ref{def_permutation}. We place these numbers in the cells of the
bottom border of the desired deco polyomino (see Figure
\ref{bij_5} b)) and for each cell with a nonzero entry $k$, we
place $k$ cells on the top of its left neighbor. The obtained deco
polyomino is the image $\Phi_5(\pi)$ of the given permutation
(Figure \ref{bij_5} c)).

The inverse map is defined in the following manner. We take a deco
polyomino of height $n$ and we place a $0$ in the first cell of
each row of the bottom border. In the remaining  cells of the
bottom border we place the number of cells situated above its left
neighbor (Figure \ref{bij_5} c)). We obtain the sequence
$b=(b_1,b_2,\ldots, b_n)$. We do have $b_i\leq n-i$ for
$i=1,2,\dots,n$. Indeed, mark the first $i$ cells along the bottom
border of the deco polyomino and note that, due to the fact that
height is attained only in the last column, the number of the
cells that define $b_i$ does not exceed the number of unmarked
cells along this bottom border. Consequently, there is a unique
permutation $\pi_1\pi_2\ldots \pi_n \in S_n$ whose right inversion
vector is $b$. If the lengths of the rows of the bottom border of
the deco polyomino are $s_1,s_2,\ldots,s_r$, then its image under
the inverse bijection, written in cycle form, is defined to be
$$
(\pi_1\pi_2\ldots \pi_{s_1})(\pi_{s_1+1}\pi_{s_1+2}\ldots
\pi_{s_1+s_2})\ldots(\pi_{s_1+s_2+\ldots+s_{r-1}+1}\ldots
\pi_{s_1+s_2+\ldots+s_{r}})
$$

\begin{figure}[!hbp]
\begin{center}
\includegraphics[scale=.38]{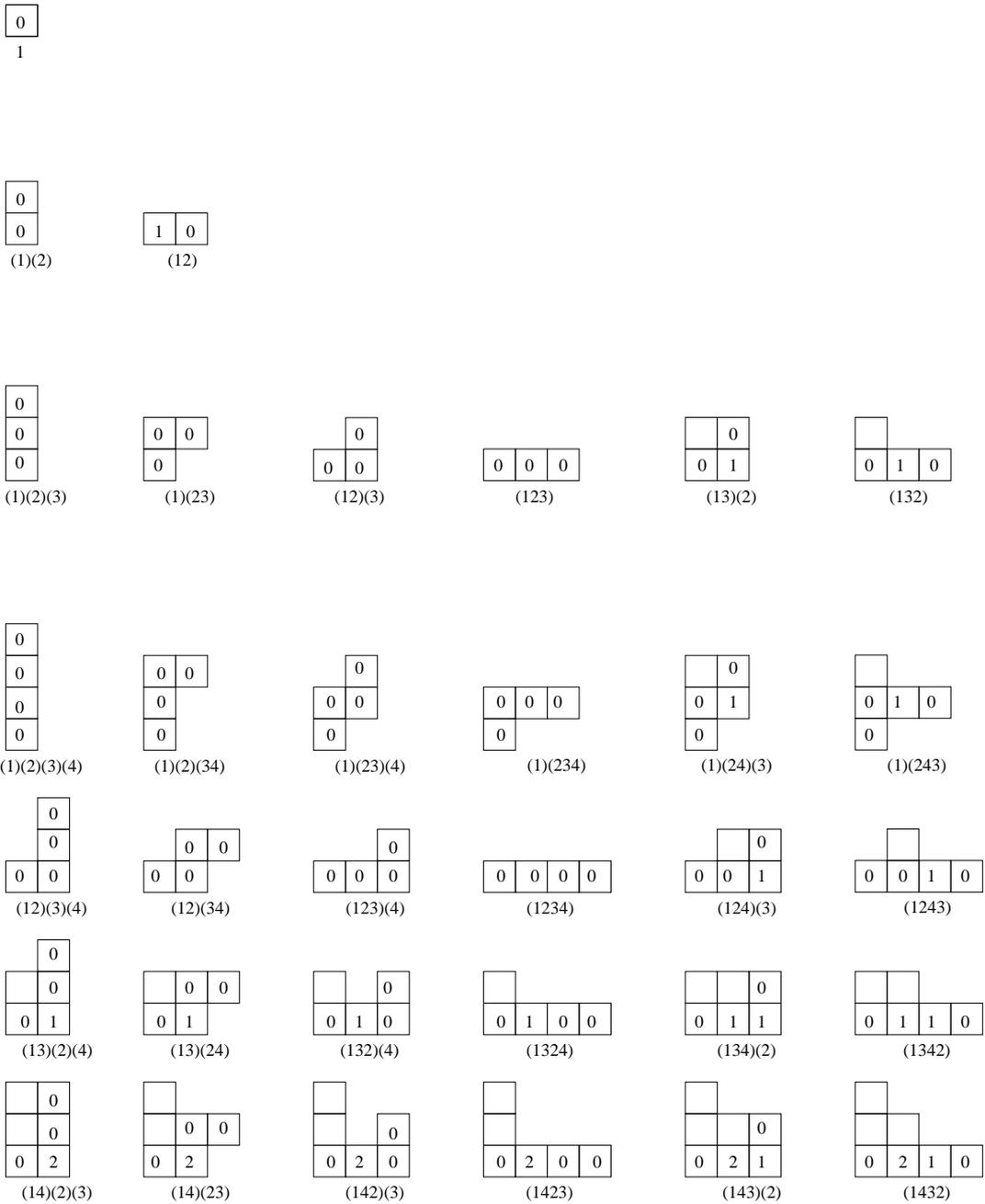}
\end{center}
\caption{Bijection $\Phi_5$ for $n=1, 2, 3, 4$} \label{ex bij_5}
\end{figure}

Figure \ref{ex bij_5} shows this bijection for $n=1, 2, 3, 4$.

\bigskip

\bigskip

\begin{example}
Consider the deco polyomino of Figure \ref{bij_5} c). Then
$b=(0,1,0,$ $3,0,0,2,1,0)$. The unique permutation with this right
inversion vector is $\pi'=132745986$. Now, taking into account
that the lengths of the rows of the bottom border of the given
deco polyomino are $5$, $3$ and $1$, we recapture the desired
permutation $\pi=(13274)(598)(6)=372196458$.
\end{example}
\bigskip

From the definition of this bijection we obtain the following
relations regarding a permutation $\pi \in S_n$ and its
corresponding deco polyomino $\Phi_5(\pi)$:
\begin{itemize}
    \item the area of $\Phi_5(\pi)$ is equal to $n+inv_c(\pi)$;
    \item the level of the last column of $\Phi_5(\pi)$ is equal to the number
    of cycles of $\pi$;
    \item the lengths of the cycles of $\pi$, in the order they are
    listed in the standard cycle form, are equal, respectively,
    to the lengths of the rows of the bottom border of $\Phi_5(\pi)$,
    starting with the lowest one.
\end{itemize}

\bigskip

\section{Bijection No. 6}\label{sec_bij_6}

\begin{figure}[!hbp]
\begin{center}
\includegraphics[scale=.4]{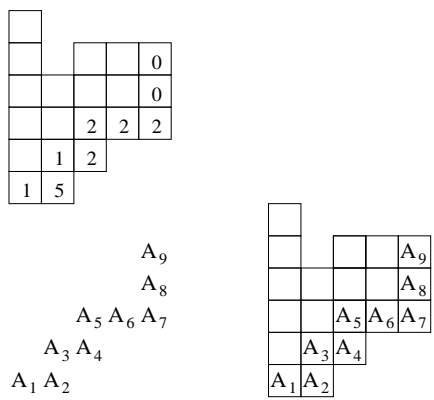}
\end{center}
\caption{An example for bijection $\Phi_6$} \label{bij_6}
\end{figure}

We define now our last bijection $\Phi_6: S_n \rightarrow D_n$. We
start with $\Phi_6^{-1}$. Parallel with its presentation for an
arbitrary deco polyomino  of height $n$, we exemplify its steps on
the polyomino of Figure \ref{bij_6}.

\begin{figure}[!hbp]
\begin{center}
\includegraphics[scale=.38]{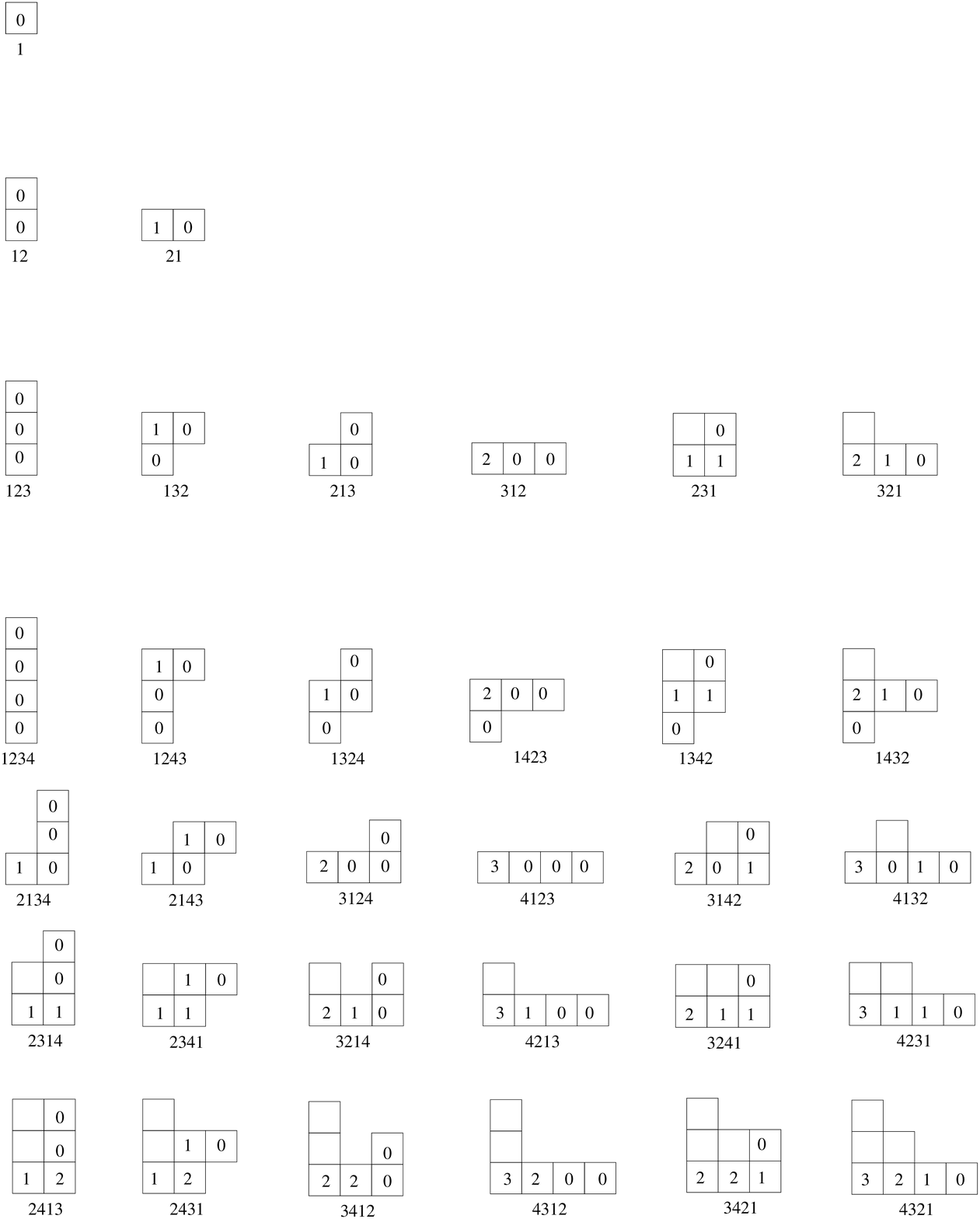}
\end{center}
\caption{Bijection $\Phi_6$ for $n=1, 2, 3, 4$} \label{ex bij_6}
\end{figure}

\bigskip
In each cell of the bottom border of the given deco polyomino we
write a number according to the following rule. In the first cell
of each row of the bottom border we write the number of cells on
its right. For the polyomino of Figure \ref{bij_6} these are the
numbers $1$, $1$, $2$, $0$ and $0$. For each of the remaining
cells of the bottom border we write the number of cells above the
cell on the left. For the deco polyomino of Figure \ref{bij_6}
these are the numbers $5$, $2$, $2$ and $2$. Next we consider the
sequence $b_1,b_2,\ldots,b_n$ of the numbers placed in the cells
of the bottom order $(1,5,1,2,2,2,2,0,0$ in the case of our
example). As in Section \ref{sec_bij_5}, one can justify that $0
\leq b_i\leq n-i$ for $i=1,\ldots,n$. Now we find the unique
permutation $\pi=\pi_1\pi_2\ldots \pi_n$ whose right inversion
vector is $(b_1,b_2,\ldots,b_n)$. This permutation $\pi$ is the
image of the $\delta$ under our bijection. For our example, we
find $\pi=273568914$.

We explain the reverse map on this example, i.e. we take

\begin{eqnarray*}
  \pi &=& 2\ 7\ 3\ 5\ 6\ 8\ 9\ 1\ 4 \\
    & & 1\ 5\ 1\ 2\ 2\ 2\ 2\ 0\ 0
\end{eqnarray*}
where below the entries $\pi_i$ we have placed the terms $b_i$ of
the right inversion vector of $\pi$. To the right of the source
cell we have $b_1=1$ cells. This allows us to determine the
positions of the first $3$ cells of the bottom border of the
desired polyomino (cells $A_1$, $A_2$, $A_3$ on Figure
\ref{bij_6}). Since $b_3=1$, there is one more cell in the second
row of the bottom border; this determines the positions of cell
$A_4$ and $A_5$ of the bottom border. Since $b_5=2$, there are 2
more cells in the third row of the bottom border, determining the
positions of the cells $A_6$, $A_7$, and $A_8$. Since $b_8=0$,
there are no other cells in the fourth row of the bottom border,
determining the position of cell $A_9$. Since $b_9=0$, there are
no cells to the right of $A_9$. Finally, we place columns of
lengths $b_2=5$, $b_4=2$, $b_6=2$, and $b_7=2$ over the cells
situated at the left of the cells $A_2$, $A_4$, $A_6$, and $A_7$,
respectively.

Figure \ref{ex bij_6} shows this bijection for $n=1, 2, 3, 4$.
\bigskip

From the definition of this bijection we obtain the following
relations regarding a permutation $\pi\in S_n$ and its
corresponding deco polyomino $\Phi_6(\pi)$:
\begin{itemize}
    \item  the number of cells of the first row of $\Phi_6(\pi)$ is equal to  $\pi_1$ (the first entry of the
    permutation);
    \item the area of $\Phi_6(\pi)$ $-$ the level of the last column
    of $\Phi_6(\pi)$ is equal to $inv(\pi)$.
\end{itemize}

\section{Acknowledgments}
The authors want to thank  E. Barcucci and L. Ferrari for their
useful suggestions.

\end{document}